\def\draft{y}
\documentclass{amsart}
\usepackage{fullpage,amssymb,epic,eepic,epsfig}

\theoremstyle{plain}

\newtheorem{theorem}{Theorem}
\newtheorem{proposition}{Proposition}[section]
\newtheorem{lemma}[proposition]{Lemma}

\newtheorem{conjecture}{Conjecture}

\theoremstyle{definition}
\newtheorem{definition}[proposition]{Definition}

\theoremstyle{remark}

\def\printname#1{
	\if\draft y
		\smash{\makebox[0pt]{\hspace{-0.5in}
			\raisebox{8pt}{\tt\tiny #1}}}
	\fi
}

\newlength{\standardunitlength}
\catcode`\@=11
\long\def\@makecaption#1#2{%
     \vskip 10pt

\setbox\@tempboxa\hbox{
       \small\sf{\bfcaptionfont #1. }\ignorespaces #2}%
     \ifdim \wd\@tempboxa >\captionwidth {%
         \rightskip=\@captionmargin\leftskip=\@captionmargin
         \unhbox\@tempboxa\par}%
       \else
         \hbox to\hsize{\hfil\box\@tempboxa\hfil}%
     \fi}
\font\bfcaptionfont=cmssbx10 scaled \magstephalf
\newdimen\@captionmargin\@captionmargin=2\parindent
\newdimen\captionwidth\captionwidth=\hsize
\catcode`\@=12

\def\BN{\mathbb N}
\def\BZ{\mathbb Z}
\def\BQ{\mathbb Q}
\def\BR{\mathbb R}
\def\BC{\mathbb C}
\def\BE{\mathbb E}

\def\A{\mathcal A}
\def\B{\mathcal B}
\def\C{\mathcal C}
\def\E{\mathcal E}
\def\D{\Delta}
\def\K{\mathcal K}
\def\cL{\mathcal L}
\def\G{\mathcal G}
\def\O{\mathcal O}
\def\T{\mathcal T}
\def\M{\mathcal M}
\def\N{\mathcal N}
\def\F{\mathcal F}
\def\V{\mathcal V}
\def\W{\mathcal W}
\def\P{\mathcal P}
\def\R{\mathcal R}
\def\T{\mathcal T}

\def\I{\mathcal I}
\def\CI{\mathcal C}
\def\Err{\mathcal Err}

\def\a{\alpha}
\def\bb{\beta}
\def\La{\Lambda}
\def\l{\lambda}
\def\Ga{\Gamma}
\def\S{\Sigma}
\def\Sp{\Sigma^p}
\def\s{\sigma}
\def\varsi{\varsigma}
\def\ga{\gamma}

\def\as{algebraically split}
\def\ihs{integral homology 3-sphere}
\def\qhs{rational homology 3-sphere}
\def\hs{homology 3-sphere}
\def\fti{finite type invariant}
\def\uf{unit-framed}
\def\rf{rationally framed}
\def\Heg{Heegaard}
\def\cok{\mathrm{cok}}
\def\torsion{\mathrm{torsion}}

\def\la{\langle}
\def\ra{\rangle}

\def\Fb#1{\mathcal F^{b}_{#1}\mathcal M}  
\def\FY#1{\mathcal F^Y_{#1}\mathcal M}

\def\Gb#1{\mathcal G^{b}_{#1}\mathcal M}  
\def\GY#1{\mathcal G^Y_{#1}\mathcal M}

\def\w{\omega}
\def\ov#1{\overline{#1}}

\def\we{\wedge}
\def\gl{\mathfrak{gl}}

\def\g{\gamma}
\def\t{\tau}
\def\e{\epsilon}
\def\Ga{\Gamma}
\def\d{\delta}
\def\b{\beta}
\def\th{\theta}
\def\Th{\Theta}
\def\s{\sigma}
\def\p{\prime}
\def\sub{\subseteq}

\def\lk{{\text{lk}}}
\def\sgn{\operatorname{sgn}}

\newcommand{\Ker}{\operatorname{Ker}}
\newcommand{\con}{\equiv}
\def\Sei{\mathrm{Sei}}

\def\iso{\cong}
\def\pt{\mathbb Z[\![t_1 ,\dots,t_m ]\!]}
\def\sms{\smallsmile}       
\def\smf{\smallfrown}       

\def\sminus{\smallsetminus}

\def\ti{\widetilde}

\def\circle{\operatorname{\psfig{figure=draws/circle.eps,height=0.1in}}}
\def\line{\operatorname{\psfig{figure=draws/line.eps,height=0.1in}}}
\def\wheel{\operatorname{\psfig{figure=draws/w2.eps,height=0.1in}}}
\def\ostar{\circledast}
\def\tcircle{\circlearrowleft}
\def\tline{\uparrow}

\def\strutb#1#2#3{\overset{#1}{\underset{#2}{ 
\begin{array}{c} \vspace{0.0cm}
\uparrow 
\vspace{-0.25cm} \\        
| \vspace{-0.45cm} \\      
\bullet \vspace{0.00cm}   
\end{array} }}\! #3}

\def\st#1#2#3{\overset{#1}{\underset{#2}{
\begin{array}{c} \hspace{-1.3mm}
	\raisebox{-4pt}{\psfig{figure=draws/strut.eps,height=0.2in} }
	\hspace{-1.9mm}\end{array} }}#3}

\def\AS{\mathrm{AS}}
\def\IHX{\mathrm{IHX}}
\def\GA{\A^{\mathrm{gp}}}
\def\GAz{\A^{\mathrm{gp},0}}

\def\mat#1#2#3#4{\left(
\begin{matrix}
 #1 & #2  \\
 #3 & #4   
\end{matrix}
\right)}

\def\matt#1#2#3#4#5#6#7#8#9{\left[
\begin{matrix}
 #1 & #2 & #3 \\
 #4 & #5 & #6 \\
 #7 & #8 & #9   
\end{matrix}
\right]}

\def\Sym{\mathrm{Sym}}
\def\pt{\partial}

\def\Per{\mathrm{Per}}
\def\fsl{\mathfrak s\mathfrak l}
\def\LIG{\mathrm{LIG}}
\def\LID{\mathrm{LID}}
\def\TID{\mathrm{TID}}
\def\IM{\mathrm{IM}}
\def\IMJ{\mathrm{IM}_J}
\def\ot{\otimes}
\newcommand{\cl}{$\clubsuit$}
\def\bt{\bar{t}}
\def\rot{\mathrm{rot}}
\def\longto{\longrightarrow}
\def\Sev{{\mathcal S}}

\def\Ges{\G^{ev,s}}
\def\sign{\mathrm{sign}}
\def\exc{\mathrm{exc}}
\def\Cab{\mathrm{Cab}}
\def\dgn{\mathrm{dgn}}
\def\udot{\mathaccent\cdot\cup}
\def\mult{\mathrm{mult}}
\def\as{\mathrm{as}}
\def\smooth{\mathrm{smooth}}
\def\bq{\bar{q}}
\def\err{\mathrm{err}}
\def\dec{\mathrm{dec}}
\def\vert{\mathrm{vert}}
\def\circ{\mathrm{circ}}
\def\longto{\longrightarrow}
\def\unknot{\mathrm{unknot}}
\def\defect{\mathrm{def}}
\def\type{\mathrm{type}}
\def\Per{\mathrm{Per}}
\def\vphi{\varphi}
\def\Edges{\mathrm{Edges}}

\begin{document}


\title[Ground State Incongruence in 2D Spin Glasses Revisited]
{Ground State Incongruence in 2D Spin Glasses Revisited}

\author{Martin Loebl}
\address{Department of Applied Mathematics \\
and \\
Institute for Theoretical Computer Science (ITI)\\
Charles University \\
Malostranske n. 25 \\
118 00 Praha 1 \\
Czech Republic.}
\email{loebl@kam.mff.cuni.cz}

\thanks{Gratefully acknowledges the support of ICM-P01-05.
This work was done while visiting the DIM, U. Chile.}


\date{
This edition: \today \hspace{0.5cm} First edition: July 20, 2003.}


\begin{abstract}
A construction supporting a conjecture that different ground state pairs 
exist in the 2-dimensional Edwards-Anderson Ising spin glass is presented.
\end{abstract}

\maketitle



\section{Introduction}
\label{sec.intro}
A fundamental and extensively studied problem on the way towards understanding the effects of disorder and frustration is to determine the multiplicity of infinite volume groundstates in 
finite-dimensional realistic models. One conjecture, in analogy with the mean-field 
Sherrington-Kirkpatrick model, is that finite-dimensional short-ranged systems with 
frustration have infinitely many groundstate pairs (\cite{MPV}, \cite{BY}). A different conjecture based on droplet-scaling theories predicts that only one groundstate pair exists (\cite{M}, \cite{BM}, \cite{FH}). 

The simplest system used to study these questions is the Edwards-Anderson Ising model (\cite{EA}) in dimension two. Here the hypothesis that  only one groundstate pair exists has received support from seminal analytic work 
 of Newman and Stein (\cite{NS2}, \cite{NS1}).  The purpose of this paper is to present a combinatorial construction supporting the competing hypothesis. 
In particular very intuitive Conjectures \ref{conj.2}, \ref{conj.b2} are
 formulated whose validity implies that incongruent (finitely incongruent respectively) groundstate pairs exist.  
These conjectures concern only finite sublattices of the square lattice and hence they may be studied by many tools including computer simulations.

The Edwards-Anderson Ising model on a graph $G=(V,E)$ 
is defined by coupling constants $J_{ij}$ assigned to each edge $\{i,j\}$ of $G$. 
We will assume that $J_{ij}$'s are independently chosen 
from a mean zero Gaussian distribution.
A physical state of the system is given by a {\em spin
assignment} $\sigma: V \rightarrow \{\pm 1\}$ which has the corresponding
energy $$ E(\sigma) = -\sum_{\{i,j\} \in E}{J_{ij} \sigma_i \sigma_j}.$$
A state is {\it groundstate} if its energy cannot be lowered by changing an arbitrary
finite set of spins. Groundstates exist for the square lattice ${\mathcal S}$ and
arbitrary coupling constants $J_{ij}$ assigned to its edges by a compactness
argument. Note that if we reverse all spins in a groundstate we again get a groundstate.
Let us call these pairs {\it groundstate pairs}, or GSPs. 
Edge $ij$ is {\it satisfied} by spin assignment $\sigma$ if $J_{ij}\sigma_i \sigma_j > 0$.
Two GSPs are called {\it incongruent} if the set of edges  satisfied by exactly one of them
 has a positive density. Note that the connectivity components of such a set in 
the dual lattice are sometimes called  'domain walls'. 

In 1D there is no frustration and only a single GSP exists. In other dimensions the main incongruency problem may be formulated as follows:

\begin{conjecture}
\label{conj.main}
Let the coupling constants $J_{ij}$ in the square lattice be chosen at random. Then with probability strictly bigger than zero there are incongruent GSP's.
\end{conjecture}

In their strategy to prove that incongruent GSP's do not exist in the 2-dimensional square
lattice, Newman and Stein (\cite{NS2}, \cite{NS1}) approach the main problem
 by means of {\it metastates}. A metastate has been introduced as a translation invariant
measure constructed as follows: for each finite sublattice $S_L$ of the square 
lattice $\S$ with {\em periodic boundary conditions} consider the joint distribution of 
coupling constants and GSPs in $S_L$. When $L$ goes to infinity, by compactness, 
there is a subsequence of $L$'s so that the joint distributions converge to translation-invariant 
(since periodic boundary conditions are imposed) joint measure. The metastate induces a 
translation invariant measure on the sets of edges  satisfied by exactly one of two GSPs. In 
this setting Newman and Stein formulate their conjecture as follows: 

\begin{conjecture}
\label{conj.NS}
Two randomly chosen GSPs (from the same metastate) are not incongruent.
 \end{conjecture}
 
Note an important fact (see Lemma~2 of \cite{NS1}) which follows from 
the translation invariance: if two randomly chosen GSPs from a metastate are distinct, 
 then with probability one they are incongruent. Newman and Stein give support to their 
conjecture in \cite{NS1}.  In particular they show that if  two GSPs chosen at random from a metastate are distinct then there is exactly one domain wall between them and it is 
a both-ways-infinite path. Newman and Stein consider this situation unlikely and express 
belief that all the metastates  are the same. 

In this paper we consider the finite sublattices with different
than periodic boundary conditions: we fix the spins along the boundary of the finite
sublattice so that maximum number of the edges of the boundary are satisfied. 
 Note an important fact: a state of minimum energy with these boundary conditions
need not be a groundstate. We will call it a {\it c-groundstate}. 
Apart of the traditional incongruency we also consider a weaker notion: 
we say that two states are {\it finitely incongruent} if  at least one of the domain walls between them contains an infinite both ways unbounded path. 
This is certainly a weaker notion of incongruency than the one of Newman and Stein
but, also in view of their results described above, the existence of finitely incongruent
states would support a conjecture that incongruent states exist.
 In our setting a weaker incongruency conjecture than Conjecture \ref{conj.main}
may be formulated as follows:

\begin{conjecture}
\label{conj.1}
Let coupling constants $J_{ij}$  in the square lattice $\S$ be chosen at random. Then with 
probability strictly bigger than zero there is a nested sequence $S_i, i=1,2,...$ 
of finite sublattices converging to the square lattice $\S$ so that 
if $e_j$ is a c-groundstate in $S_{2j}$ and $o_j$ is a c-groundstate in $S_{2j-1}$ 
then both sequences $(o_j)$, $(e_j)$ converge and their limits $o,e$ are finitely 
incongruent states of  $\S$.
\end{conjecture}

Note that $o,e$ are groundstates of $\S$ since the boundary conditions 
'disappear to the infinity'. Also note that the boundary conditions mean that all 
or all but one boundary edges are to be satisfied depending on the parity 
 of the number of edges with negative coupling constants.

We start by considering the strip lattice ${\mathcal C}_k$: the vertical coordinates of its vertices are arbitrary integers, and its horisontal coordinates run through 
integers from $-k$ to $k$. As an introduction to our method we show in section \ref{sec.strip} that the strip lattice satisfies Conjecture \ref{conj.1}.  In section \ref{sec.fin} we 
formulate two conjectures which imply Conjectures \ref{conj.1} and \ref{conj.main}.
The important feature of these new conjectures is that they concern finite sublattices
only, and hence allow direct study by discrete methods and simulations.
In section \ref{sec.pin} we prove a (rather weak) consequence of Conjecture \ref{conj.2}.
The discrete aspect of Conjectures \ref{conj.2}, \ref{conj.b2} is supported in the last section 
\ref{sec.dual} where we show that the dual 
formulations of the two conjectures are statements about $T$-joins in finite square lattices; 
$T$-join belongs to basic discrete optimization notions and as such it is heavily studied by 
discrete and computational methods.

\section{The strip Lattice}
\label{sec.strip}

Let $C(n,k)$ be the finite induced subgraph of ${\mathcal C}_k$  with the vertices $(i,j): |j|\leq n$. The basic building blocks of strip and square lattices are  unit squares called {\it plaquettes}. 
A plaquette is {\it frustrated} if it has odd number of edges (out of 4) with negative coupling 
constants. Observe that a plaquette is frustrated  if and only if arbitrary state satisfies an odd 
number of its edges. 
We define graph $C(n,k)^*$ whose vertices are all the plaquettes 
of $C(n,k)$ and the edges are all edges $e^*$ such that $e$ is an edge of $C(n,k)$
not on the boundary; edge $e^*$ connects two plaquettes $p,q$ such that edge $e$ lies 
on the boundary of $p$ and $q$. Note that this differs from a standard definition of the dual
planar graph since we do not consider the duals of the boundary edges.
If $A$ is a subset of edges of $C(n,k)$ then let $A^*$ denote the set of 'dual' edges: $A^*=\{e^*; e\in A\}$. Let $n>m$, and consider the graph 
$D(n,m,k)=C(n,k)-C(m-1,k)$. Note that $D(n,m,k)$ has two connectivity components, 
each of them consists of $n-m$ horizontal levels of plaquettes. The two components of $D(n,m,k)$ are naturally called {\it upper} and {\it lower} and denoted by $DU(n,m,k)$ 
and $DL(n,m,k)$. 
 
\begin{definition}
\label{def.reg}
We say that $C(n,k), C(m,k)$ is a {\it regular pair} if  for each $k'\leq k$,
 both $C(n, k')$ and $C(m, k')$ have an even number of boundary edges with negative coupling constants and both $DU(n,m,k)$ and $DL(n,m,k)$ have exactly one frustrated plaquette, located in the middle of the lowest (highest respectively) horizontal row. 
\end{definition}

We will use the following key observation:

\begin{lemma}
\label{lem.1}
Let $C(n,k), C(m,k)$ be a {\it regular pair}.
Let $c(i)$ be a c-groundstate of $C(i,k)$ and let $DIS(c(i))$ denote the set of edges 
dissatisfied by $c(i)$, $i=n,m$.Then the symmetric difference $DIS(c(n))^* \Delta DIS(c(m))^*$ contains a path from a plaquette of $DU(n,m,k)$ to a plaquette of $DL(n,m,k)$.
\end{lemma}
\begin{proof}
The subgraph formed by $DIS(c(i))^*$ induces odd degree in each frustrated plaquette
and even degree in each happy plaquette of $C(i,k)$. Moreover for $i=m,n$, $DIS(c(i))$
contains no edge of the boundary of $C(i,k)$.
Hence $DIS(c(n))^*\Delta DIS(c(m))^*$ induces odd degree 
in each frustrated plaquette of $D(n,m,k)$, and even degree in arbitrary other 
plaquette of $C(n,k)$. This easily implies the Lemma. 
\end{proof}

Now we are ready to show that 

{\bf Conjecture \ref{conj.1} holds for the strip lattice:}

Let $S_i=C(i,k)$. Clearly, for almost all coupling constants assignments $J$
in the whole strip lattice there is a sequence $(m_j)$ so that for each $j$,
$C(m_j, k), C(m_j-1, k)$ is a regular pair. 
Let $o_j$ be a c-groundstate in $S_{m_j}$ and let $e_j$ be a c-groundstate in $S_{m_j-1}$.
Lemma \ref{lem.1} implies that for each $j$, $DIS(o_j)^*\Delta DIS(e_j)^*$ contains a path $P_j$ of length at least $2j+1$. Now it is easy to see that

{\it Claim~1.} There is a subsequence $P_{n_j}$ that converges to both ways infinite
path $P$.

By compactness there is a subsequence $(p_j)$ of $(n_j)$ so that both sequences 
 $(e_{p_j})$ and $(o_{p_j})$ converge. Let the respective limits be $e$ and $o$. 
Then necessarily $P$ is a subset of a domain wall between $e$ and $o$ and so 
$e,o$ are finitely incongruent. Hence Conjecture \ref{conj.1} holds for the strip lattice.

\section{The Finite Conjectures}
\label{sec.fin}

In this section we formulate Conjecture \ref{conj.2} and Conjecture \ref{conj.b2}.
We show that Conjecture \ref{conj.2} implies Conjecture \ref{conj.1} and Conjecture \ref{conj.b2} 
implies the main Conjecture \ref{conj.main}.

\begin{definition}
\label{def.A}
We first introduce some notation.
\begin{itemize}
\item
We denote by ${\bf  A(n,k)}$ the set of all coupling constants assignments in $\S$ so that there is a c-groundstate $r$ in $C(n-1,k)$ and a c-groundstate $s$ in $C(n,k)$ and a path $P$ in $DIS(r)^*\Delta DIS(s)^*$ from a frustrated plaquette of $DU(n,n-1,k)$ to a frustrated plaquette of $DL(n,n-1,k)$ which contains an edge in distance at most $100$ from the origin.
\item
A row $R$ of plaquettes of $C(n,k)$ is called {\it isolation} if the middle plaquette $M$ is the only frustrated one and for each horizontal edge $e\notin M$, $|J_e|> \sum_{e'}|J_{e'}|$
where the sum is over all edges $e'$ such that $e'\in M$ or $e'$ is a vertical edge 
incident to a vertex of a plaquette of $R$.  $C(n,k)$ is called {\it isolated}
if both boundary horizontal levels of plaquettes are isolation.
\item
We denote by ${\bf R(n,k)}$ the set of all coupling constants assignments in $\S$ so that 
$C(n,k)$, $C(n-1,k)$ is a regular pair and $C(n,k)$ is isolated.
\item
Finally we denote by ${\bf R'(n,k)}$ the set of all coupling constants assignments in $\S$ 
which belong to $R(n,k)$ and do not belong to  $R(n',k)$ for $k\leq n'<n$.
\end{itemize}
\end{definition}

We will show that the following conjecture implies Conjecture \ref{conj.1}.   

\begin{conjecture}
\label{conj.2}
There is positive integer constant $c$ and $\epsilon> 0$ so that if $n, k> c+1$ 
then the probability of $A(n,k)$, in the set of the coupling constants assignments $J$
of $\S$ such that $J\in A(n,k')$ for each $c+1<k'<k$ and $J\in R'(n,k)$, is at least 
$1- (k-c)^{-1}(\log (k-c)^{-1-\epsilon}$.
\end{conjecture}

How is it possible that for each $k'< k$, $A(n,k')$ holds and $A(n,k)$ does not hold? 
There may be a block of heavy edges of width $2k-1$ encircling the origin.
The domain walls in $C(n,k')$, $k'<k$, pass through it since they cannot escape elsewhere
however the domain wall in $C(n,k)$ avoids it. If this is essentially the only possibility,
Conjecture \ref{conj.2} should be true.

If Conjecture \ref{conj.1} holds then it is natural to expect that for large $k$, a neighbourhood of the origin behaves in a similar way as the origin itself. This leads to a bolder Conjecture
\ref{conj.b2} which we will show implies the main Conjecture \ref{conj.main}.

\begin{definition}
\label{def.BA}
We denote by ${\bf  BA(n,k)}$ the set of all coupling constants assignments in $\S$ so that there is a c-groundstate $r$ in $C(n-1,k)$ and a c-groundstate $s$ in $C(n,k)$ and a path $P$ in $DIS(r)^*\Delta DIS(s)^*$ from a frustrated plaquette of $DU(n,n-1,k)$ to a frustrated plaquette of $DL(n,n-1,k)$ so that $P$ contains an edge in distance at most $100$
from the origin AND has density at least $1/100$ in the square centered at the origin 
with the side-length $100k^{1/100}$.
\end{definition}

\begin{conjecture}
\label{conj.b2}
Conjecture \ref{conj.2} holds also with $A(n,k')$ replaced by $BA(n,k')$, $k'\leq k$.
\end{conjecture}

\begin{theorem}
\label{thm.1}
  Conjecture \ref{conj.2} implies  Conjecture \ref{conj.1}.
\end{theorem}
\begin{proof}

We prove this theorem in a series of observations. 

{\it Observation~1:} Let $C(n,k) \in R(n,k)$.
 Let $r$ be a c-groundstate in $C(n-1,k)$, $s$ be a c-groundstate  in $C(n,k)$ and let $P$
be a path in $DIS(r)^*\Delta DIS(s)^*$ from the frustrated plaquette of $DU(n,n-1,k)$ to  the frustrated plaquette of $DL(n,n-1,k)$. Moreover let $n'> n$, let $C(n',k)$, $C(n'-1,k)$ be a regular pair and let $r', s', P'$ be defined analogously as $r, s, P$. Then $P'$ contains $P$.

Let $k> c+ 1$ and let $I_k$ denote the set of all coupling constants assignments in $\S$ so that for each $c+1< k'\leq k$ there is $n(k')\geq k'$ so that  $C(n(k'),k')  \in R(n(k'),k') \cap 
A(n(k'),k')$. 

{\it Observation~2:} If $k> c+ 1$ and Conjecture \ref{conj.2} holds then 
$$
Prob(I_k)= Prob(\cap_{c+1< k'\leq k}I_{k'})\geq 
\prod_{j= c+2}^k(1- (j-c)^{-1}(\log (j-c))^{-1-\epsilon}).
$$
{\it Proof.}

 First notice that for each $k$ with probability one there is an infinite sequence $(n_i)$ such that $R(n_i, k)$ holds for each $i$. In particular if we let $Z(k)$ be the event 
'There is no $n\geq k$ with $R(n, k)$', then $Z(k)$ has probability zero.

We proceed by induction on $k$. The case $k= c+2$ follows from the fact above and Conjecture \ref{conj.2}.  For the induction step first note that event $Z(k)$ has probability zero and so it also has probability zero conditioned on $I_{k-1}$, if we use the induction assumption. Hence for almost all elements of $I_{k-1}$ there is smallest $n\geq k$ such that
$R(n,k)$ holds. Next note that Observation~1 implies that the set of instances satisfying 
$R'(n, k)$ and $I_{k-1}$ is the same as the set of instances satisfying $R'(n, k)$ and 
$A(n,k')$ for each $c+1< k'< k$.  Hence assuming validity of Conjecture \ref{conj.2}
and the induction assumption we get
$$
Prob(I_k)= Prob(I_{k-1}) Prob(I_k || I_{k-1})=
Prob(I_{k-1}) \sum_m Prob(R'(m,k)|| I_{k-1})Prob(I_k ||R'(m,k), I_{k-1})=
$$
$$
Prob(I_{k-1}) \sum_m Prob(R'(m,k)|| I_{k-1})Prob(I_k ||R'(m,k), A(m,k'), c+1<k'<k)=
$$
$$
Prob(I_{k-1}) \sum_m Prob(R'(m,k)|| I_{k-1})Prob(A(m,k) ||R'(m,k), A(m,k'), c+1<k'<k)\geq
$$
$$ 
\sum_mProb(R'(m,k)|| I_{k-1})\prod_{j= c+2}^k(1- (j-c)^{-1}(\log (j-c))^{-1-\epsilon})=
$$
$$
\prod_{j= c+2}^k(1- (j-c)^{-1}(\log (j-c))^{-1-\epsilon}).
 $$
 This finishes the proof of Observation~2.

Hence assuming validity of Conjecture \ref{conj.2} the probability of the event
 'For each $k$, $I_k$' is at least $\prod_{j\geq 2}(1- j^{-1}(\log j)^{-1-\epsilon}> 0$. This 
proves Theorem \ref{thm.1} in the same way as Claim~2 proves Conjecture \ref{conj.1}
for the strip lattice.
\end{proof}

\begin{theorem}
\label{thm.b1}
  Conjecture \ref{conj.b2} implies  Conjecture \ref{conj.main}.
\end{theorem}
\begin{proof}
We proceed similarly as in the proof of the previous theorem.
Let $k> c+ 1$ and let $BI_k$ denote the set of all coupling constants assignments in $\S$ so that for each $c+1< k'\leq k$ there is $n(k')\geq k'$ so that  $C(n(k'),k')  \in R(n(k'),k') \cap 
BA(n(k'),k')$.  We observe as before if $k> c+ 1$ and Conjecture \ref{conj.b2} holds then 
$$
Prob(BI_k)= Prob(\cap_{c+1< k'\leq k}BI_{k'})\geq 
\prod_{j= c+2}^k(1- (j-c)^{-1}(\log (j-c))^{-1-\epsilon}).
$$
Hence assuming validity of Conjecture \ref{conj.b2} the probability of the event
 'For each $k$, $BI_k$' is at least $\prod_{j\geq 2}(1- j^{-1}(\log j)^{-1-\epsilon}> 0$. This 
proves Theorem \ref{thm.b1}: in the same way as above we can grow a path
in the symmetric difference, so that it goes near to origin AND gradually 'fills'
the whole square grid.
\end{proof}

\section{The Pinning Lemma}
\label{sec.pin}

In the rest of the paper we collect support for Conjecture \ref{conj.2}. First we present the following Pinning Lemma.

\begin{lemma}
\label{lem.pin1}
There is a function $g$ from positive integers to $(0,1)$ so that for each $k$ and $n>k$,
if coupling constants in $C(n,k)$ are chosen at random so that $C(n,k), C(n-1,k)$
is a regular pair then the probability of $A(n,k)$ is at least $g(k)$.
\end{lemma}
We postpone the proof to the appendix. Anyway, the present proof is not satisfactory.  Although $g$ may well be a constant, at present we are able to show only a very weak
inverse exponential lower bound for it. Next we prove a consequence of Conjecture \ref{conj.1}, using the Pinning Lemma.

\begin{theorem}
\label{thm.2}
Let $(S_i= C(l_i,p_i), l_i > p_i)$ be a nested sequence of finite sublattices monotonically 
converging to the square lattice and such that for each $k$ there are sufficiently
many lattices with width $k$. Let ${\mathcal J}_i$ be the distribution 
of the coupling constants in $S_i$. Then for almost all $(J_i)_{i\geq 1}$ from 
$({\mathcal J}_i)_{i\geq 1}$ {\it there is a converging subsequence} $(J_{m_j})$
 with the following property:
if $e_j$ is a c-groundstate in $S_{m_j}$ and $o_j$ is a c-groundstate in $S_{m_j-1}$ 
(with coupling constants given by $J_{m_j}$) then both sequences $(o_j)$, $(e_j)$ converge 
and their limits $o,e$ are weakly incongruent states. 
\end{theorem}

 Theorem \ref{thm.2} follows from Claim~2 below in the same way as Conjecture 
\ref{conj.1} for the strip lattices follows from Claim~1.

{\bf Claim~2.}
For almost all $(J_i)_{i\geq 1}$ from $({\mathcal J}_i)_{i\geq 1}$ there is a subsequence 
$(C(n_i,k_i))$  of $(S_i)$ so that for each $i$, $C(n_i,k_i)$ and $C(n_i-1,k_i)$ 
with the coupling constants given by $J_{n_i}$ is a regular pair and there is a path $P_i$ 
in $DIS(o_i)^* \Delta DIS(e_i)^*$ from a frustrated plaquette of $DU(n_i,n_i-1,k)$ to
a frustrated plaquette of $DL(n_i,n_i-1,k)$ which contains an edge in distance at most $100$ 
from the origin. Here $o_i$ is a c-groundstate in $C(n_i-1,k_i)$ and $e_i$ is a c-groundstate in $C(n_i,k_i)$.

\begin{proof}
 Let $k$ be an arbitrary positive integer. Since $p_i=k$ for sufficiently many $i$'s, 
we know by the Pinning Lemma that the probability that one of $C(l_i,k)$ satisfies 
the properties of Claim~2 is very large. Hence the set of instances $(J_i)_{i\geq 1}$ 
from $({\mathcal J}_i)_{i\geq 1}$ for which the propertiess of Claim~2 donot hold has measure zero. This finishes the proof of Claim~2 and Theorem~2.
\end{proof}

Claim~2 is a consequence of the Pinning Lemma and {\it the fact that in each $S_i$ 
we assign the coupling constants independently}. The remaining obstacle in proving the full
 Conjecture \ref{conj.1} is that because of dependancies the Pinning Lemma cannot be used independently in each $S_i$. Conjecture \ref{conj.2} may be viewed as an attempt to make 
the dependencies work for us. 

\section{The Dual Formulation}
\label{sec.dual}

It seems very natural to formulate the Pinning Lemma as a property of the dual lattice. In doing so we connect our considerations with the concept of $T$-joins extensively studied in discrete optimization. This may be particularly useful for studying Conjectures \ref{conj.2},
\ref{conj.b2} computationally. In fact, the Pinning Lemma is proved in its dual form in the 
appendix. We start by listing some simple properties of lattices $C(n,k)$.
\begin{enumerate}
\item [1.]
$C(n,k)$ has an even number of negative coupling constants on the boundary if and only if
it has an even number of frustrated plaquettes.
\item [2.]
A set $R$ of edges of $C(n,k)$ not on the boundary is the set $DIS(r)$ 
of the dissatisfied edges of a state $r$ (not necessarily a groundstate) if and only if 
$R$ has an odd number of edges from each frustrated plaquette and an even number 
of edges from any other plaquette.
\item [3.]
A state $r$ is a c-groundstate if and only if it satisfies the boundary conditions and
 $\sum_{(ij)\in DIS(r)}|J_{ij}|$ is as small as possible. Hence
there is a natural bijection between the c-groundstate pairs of $C(n,k)$ and the sets 
$A$ of edges not on the boundary and satisfying: a plaquette has an odd number of edges 
of $A$ if and only if it is frustrated, and $\sum_{(ij)\in A} |J_{ij}|$ is as small as possible.
\end{enumerate}

This means that regarding the Pinning Lemma we need only a subset of information given 
by the coupling constants: we need to know the value $|J_{ij}|$ for each edge $(ij)$ 
not on the boundary, and we need to know which plaquette is frustrated. Each plaquette 
is equally likely to be frustrated or happy. If ${\mathcal J}$ is our distribution of coupling 
constants then we denote by $|{\mathcal J}|$ the distribution of their absolute values. 
We are interested only in those $C(n,k)$ that contain an even number of frustrated plaquettes. 
Hence instead of choosing the coupling constants from ${\mathcal J}$, we can choose them
 from $|{\mathcal J}|$ and choose uniformly at random an even set of plaquettes which we 
want to be frustrated. This means that the Pinning Lemma is about $C(n,k)^*$ rather than about $C(n,k)$.  $C(n,k)^*$ is also a square grid, of width $2k$ and height $2n$.
 We need one more definition. 

\begin{definition}
\label{def.T}
Let $G=(V,E)$ be a graph and let $T$ be a subset
of an even number of vertices of $G$. We say that a set $A$ of edges of $G$
is a $T$-join if each vertex $x$ of $G$ is incident with an even number 
of edges of $A$ if and only if $x\notin T$.
\end{definition}

 Taking these considerations into account, note that 
the following  Lemma \ref{lem.pin2} is an equivalent dual formulation of Pinning Lemma \ref{lem.pin1}, and Conjectures \ref{conj.3}, \ref{conj.b3} are equivalent dual formulations 
of Conjectures \ref{conj.2}. \ref{conj.b2}.

We denote by $C'(n,k)$ the graph obtained from $C(n,k)$ by attaching two more vertices
$[0,n+1]$ and $[0,-n-1]$ by the corresponding two vertical edges. We say that
$C(n,k)$ is dually isolated if each of its two boundary rows of vertices are dual isolation.
A row $R$ of vertices is dual isolation if it has exactly one vertex $r$ of $T$, located
in the middle of $R$, and the weight of each vertical edge incident to a vertex of $R-r$
is bigger than the sum of the weights of the horizontal edges in $R$ and in the two rows adjacent to $R$,
and the edges adjacent to $r$.

\begin{definition}
\label{def.D}
Let $D(n,k)$ denote the following property: there is a minimum $T$-join 
$r$ in $C(n,k)$ and a minimum $(T\cup \{[0,n+1], [0,-n-1]\})$-join $s$ in $C'(n,k)$ 
so that a 
path $P$ from $[0,n+1]$ to $[0,-n-1]$ in $r\Delta s$ contains an edge in distance at most $100$ from the origin.
\end{definition}

\begin{definition}
\label{def.BD}
Let $BD(n,k)$ denote the following property: there is a minimum $T$-join 
$r$ in $C(n,k)$ and a minimum $(T\cup \{[0,n+1], [0,-n-1]\})$-join $s$ in $C'(n,k)$ so that a 
path $P$ from $[0,n+1]$ to $[0,-n-1]$ in $r\Delta s$ contains an edge in distance at most $100$ 
from the origin AND has density at least $1/100$ in the square centered at the origin 
with the side-length $100k^{1/100}$.
\end{definition}

\begin{conjecture}
\label{conj.3}
There is positive integer constant $c$ and $\epsilon> 0$ so that if $n, k> c+1$,
and the absolute values of the coupling constants in $C(n,k)$ 
are chosen at random from $|{\mathcal J}|$, and  a subset $T$ of vertices of 
$C(n,k)$ be chosen uniformly at random so that:
 \begin{itemize}
\item
for each $k'\leq k$, $|T\cap C(n,k')|$ is even, 
\item
if $n'\leq n$ then $C(n',k)$ is not dually isolated,
\item
for each $c+1< k'< k$, $D(n,k')$ holds,

\end{itemize}
then the probability of $D(n,k)$ is at least  $1- (k-c)^{-1}(\log (k-c))^{-1-\epsilon}$.
\end{conjecture}

\begin{conjecture}
\label{conj.b3}
Conjecture \ref{conj.3} holds with $D(n,k')$ replaced by $BD(n,k')$, $k'\leq k$.
\end{conjecture}

\begin{lemma}
\label{lem.pin2}
There is a function $g$ from positive integers to $(0,1)$ so that for each $k$ and $n>k$,
if we choose absolute values of coupling constants of $C(n,k)$ at random from 
$|{\mathcal J}|$ and choose a subset $T$ of vertices of $C(n,k)$ uniformly
 at random so that for each $k'\leq k$, $|T\cap C(n,k')|$ is even then the probability 
of $D(n,k)$ is at least $g(k)$.
\end{lemma}

{\it Conclusion.}
In this paper we formulate Conjectures \ref{conj.2}, \ref{conj.b2} whose validity implies that finitely incongruent 
and incongruent groundstate pairs exist in the 2-dimensional Edwards-Anderson Ising spin glass. 
The conjectures deal with finite sublattices only and may be naturally and 
effectively  studied by computer simulations.
We gather supporting evidence, namely we prove a Pinning Lemma and as its
consequence a weaker statement in Theorem \ref{thm.2}. 

\appendix

\section{ Proof of the Dual Pinning Lemma \ref{lem.pin2}}

We will consider set ${\mathcal K}$ of {\it configurations} with joint distribution ${\mathcal U}_k$. 
A configuration is a quadruple $(J,T,x,y)$ where $J$ 
consists of the coupling constants, $T$ is an even subset of vertices of $C(n-1,k)$
and $x=[i,n], y=[j,-n]$ is a particular choice of sets $X,Y$. 
We show that there is a function $F$ from ${\mathcal K}$ to itself such that
\begin{itemize}
\item 
${\mathcal U}_k(F({\mathcal K}))\geq c_k>0$, and

\item 
Each $L\in F({\mathcal K})$ is {\it positive}, i.e. there is a minimum $T$-join $r$ in $C(n-1,k)$ 
and a minimum $T\cup \{x,y\}$-join $s$ in $C(n,k)$ such that path $P$ from $x$ to $y$ 
in $r\Delta s$ contains an edge in distance at most $100$ from the origin. 
\end{itemize}

Fix an arbitrary positive configuration $K_0$. Let $K=(J,T,x,y)$ be a configuration.
If there is an edge $e$ incident to a vertex $[i,j]$ with $|i|\leq 4$, $|j|\leq k$ 
and $|J_e|\geq 1$ then let $f(K)=K_0$, otherwise let $P, r, s$ be as in the statement of Lemma~2.
If $P$ passes in distance at most 100 from the origin then let $f(K)=K$.
Otherwise let $m$ be the smallest positive integer such that $P$ contains a vertical edge
with the x-coordinate of its vertices equal to $m$ or $-m$ and with the absolute value
of both y-coordinates at most 4. Note that $P$ has no vertex $[i,j]$
with $|i|<m$ and $|j|<4$. Let $Z$ be the graph induced on the vertices $[i,j]; |i|\leq m, |j|\leq 4$.
We let $K=K_1$, $P=P_1$, $r=r_1$, $s=s_1$, $Z=Z_1$, $T=T_1$, $m=m_1$, let $n_1$ be the number
 of vertical edges of $P\cap Z$ and let $p_1$ be the number of (all) edges of $P\cap Z$. 

Next we describe a procedure with input $I_i=(K_i=(J_i,T_i,x,y),P_i,r_i,s_i,Z_i,m_i,p_i)$
which produces $F(K)$ or $I_{i+1}$.

{\it The Procedure.} If $r_i\Delta s_i$ contains a cycle then let $F(K)=K_0$. 
Otherwise let $H_i$ be the segment of $P_i\cap Z_i$ defined as follows:
If $P_i$ contains a vertical edge $e=\{[z,a],[z,b]\}$ so that $|z|=m_i$ and $|a|<|b|<4$
 then let $H_i$ consist of $e$. If $P_i$ contains no such vertical edge but it does contain a horizontal 
edge $e=\{[a,z],[b,z]\}$ so that $|z|=4$ and $|a|<|b|<m_i$ then again let $H_i$ consist of $e$.
Finally let there be only 'corner' edges in $P_i\cap Z_i$. Let $e$ be such vertical edge
(it exists by the choice of $m_i$), $e=\{[z,a],[z,b]\}$, $|z|=m_i$ and without loss of generality
$z=-m_i, a=3, b=4$. Then we let $H_i$ consist of $e$ if $\{[z,4],[z-1,4]\}\notin P_i$, 
and $e$ together with $\{[z,4],[z-1,4]\}$ otherwise.
Let $W_i$ be the set of edges of a path in $Z_i$ between the end-points of $H_i$ 
such that it contains some edges in distance at most 100 from the origin, no vertical edge 
of the boundary of $Z_i$, and as few horizontal edges of the boundary of $Z_i$ as possible 
(i.e. at most two).
Let $M_i$ be an integer upper bound of the coupling constants of the edges incident with a vertex 
of $Z_i$. For instance $M_1=1$. For each edge $e$ of $Z_i$ such that $e\notin W_i\cup (P_i-H_i)$
we let $(J_{i+1})_e= (J_i)_e+100 k M_i$, and we let $(J_{i+1})_e= (J_i)_e$ otherwise.

$T_{i+1}$ is defined as folows: let $r'_i$ be obtained from $r_i$ by deleting all the edges of $r_i\cap s_i$
which belong to $Z_i$. Analogously define $s'_i$. Let $U^b_i$ ($U^0_i$ respectively) be the set 
of vertices of $Z_i-T_i$ ($Z_i\cap T_i$ respectively) such that we deleted an odd number 
of edges of $r_i$ incident with them. 
 We let $T'_{i+1}=(T_i-U^o_i) \cup U^b_i$. Observe that $T'_{i+1}$ has no vertices
in the interior of $Z_i$ and $r'_i$ is a $T'_{i+1}$-join and $s'_i$ is a $T'_{i+1}\cup \{x,y\}$-join.
If $H_i\subset r'_i$ or $H_i\subset s'_i$ then let $T_{i+1}=T'_{i+1}$ else necessarily $H_i$ contains two
edges incident to a 'corner vertex' of $Z_i$ and without loss of generality assume that the vertical edge
of $H_i$ belongs to $r'_i$. In this case $T_{i+1}$ is obtained from $T'_{i+1}$
by changing the status of both vertices of the horisontal edge of $H_i$. We also modify $r'_i$ and $s'_i$
so that we delete the horisontal edge of $H_i$ from $s'_i$ and add it to $r'_i$. 
Observe that  $T_{i+1}$ has no vertices in the interior of $Z_i$ and $r'_i$ is a $T_{i+1}$-join 
and $s'_i$ is a $T_{i+1}\cup \{x,y\}$-join. Moreover $H_i\subset r'_i$ or $H_i\subset s'_i$,
$r'_i\Delta s'_i=P_i$ and $r'_i\cap s'_i$ has no edges in $Z_i$.
Without loss of generality assume $H_i\subset r'_i$. Let $r''_i$ be obtained from $r'_i$ by exchanging
$H_i$ for $W_i$. Clearly $r''_i$ is a $T_{i+1}$-join, $r''_i\Delta s'_i$ is a path $P'_i$ 
obtained from $P_i$ by exchanging $H_i$ for $W_i$ and 
$E(s'_i)\leq E(s_i)$ and $E(r''_i)\leq E(r_i)+ 32M_ik$. Let $K_{i+1}=(J_{i+1},T_{i+1},x,y)$, 
$r_{i+1}$ be a minimum $T_{i+1}$-join, $s_{i+1}$ be a minimum $T_{i+1}\cup \{x,y\}$-join, and
let $P_{i+1}$ be the $x,y$-path in the symmetric difference of $r_{i+1},s_{i+1}$.
If $P_{i+1}$ contains an edge in distance at most 100 from the origin then let $F(K)=K_{i+1}$ 
otherwise we output 
vector $(K_{i+1}=(J_{i+1},T_{i+1},x,y),P_{i+1},r_{i+1},s_{i+1},S_{i+1},m_{i+1},n_{i+1}, p_{i+1})$.
This finishes the describtion of the Procedure.

Now observe that $R=r_{i+1}\Delta (r_i\Delta r'_i)$ is a $T_i$-join such that $R$ and 
$r_{i+1}$ differ only on the edges of $Z_i$. Hence
$r_{i+1}$ cannot contain an edge of $Z_i-(W_i\cup (P_i-H_i))$ since
otherwise $E(r_{i+1})\geq E(R)- 32kM_i + 100kM_i \geq E(R) + 60kM_i \geq E(r_i) + 60kM_i > E(r''_i)$. 
Hence $r_{i+1}$ contains all the edges of $W_i$ of the interior of $Z_i$ or none of them, 
and its edges from the boundary of $Z_i$ form a subset  of $Z_i\cap [(P_i-H_i)\cup D]$, where
$D$ is non-empty only if $H_i$ has a vertical 'corner' edge and then $D$ consists of
one or two horizontal edges (by the definition of $W_i$). The same holds for $s_{i+1}$.
If $H_i$ has at least one vertical edge of $Z_i$ then we have that $m_i<m_{i+1}$ or 
$m_i = m_{i+1}, n_i> n_{i+1}$. If $H_i$ has no vertical
edge then it consists of exactly one horizontal edge and all the edges of $W_i$ belong to the
interior of $Z_i$. Hence $r_{i+1}$ contains all the edges of $W_i$ or none of them
and the edges of $r_{i+1}$ from the boundary of $Z_i$ belong of $P_i-H_i$. 
The same holds for $s_{i+1}$. Hence we have $m_i<m_{i+1}$ or 
$m_i = m_{i+1}, n_i = n_{i+1}, p_i> p_{i+1}$. Sumarising, $m_i<m_{i+1}$ or 
$m_i = m_{i+1}, n_i> n_{i+1}$, or $m_i = m_{i+1}, n_i = n_{i+1}, p_i> p_{i+1}$. 

Hence after at most $(16k)^2$ repetitions of the Procedure we have $F(K)$ defined. 
Moreover $F(K)$ is defined
only if the path with desired properties exists. Finally $c_k>0$ exists since the set ${\mathcal Z}$
of configurations with $|J_{ij}|<1$ for $|i|\leq 4, |j|\leq k$ clearly satisfies 
${\mathcal U}_k({\mathcal Z})\geq c'_k>0$ and with probability one $L=F(K)$ from 
${\mathcal Z}$ may be viewed as
$K+\alpha_K$, where each component of $\alpha_K$ is a bounded integer and the number of 
possible $\alpha_K$'s is bounded from above by the number of paths on vertices 
$[i,j], |i|\leq k, |j|\leq 4$ (which is a modest function of $k$). This finishes the proof.

\vspace{10pt}
{\bf Acknowledgement.}
This project has been supported by Project LN00A056 and GAUK 158.
I have started to work on it in the autumn of 2000 while I was
 visiting Bruno Nachtergaele in Davis: I would like to thank to him for 
introduction to this subject. I am indebt to Jirka Matousek for many discussions
and in particular for suggesting a proof method for Pinning Lemma~2.
I would also like to thank to Greg Kuperberg, Jan Vondrak, Michael Lacey, Prasad Tetali, 
Russell Lyons, Laszo Erdes, Daniel Stein, Charles Newman and Matteo Palassini
 for helpful discussions.

\end{document}